\documentclass[a4 paper, 10pt]{article}
\usepackage{amsmath}
\usepackage{latexsym}
\usepackage{amsfonts}
\usepackage{amssymb}
\usepackage{amsthm}
\usepackage{amsbsy}
\usepackage{graphicx}
\usepackage{bbm}
\usepackage{placeins }
\usepackage{color}
\usepackage{epstopdf}

\newcommand{\R}{\mathbb{R}}

\newcommand{\Z}{\mathbb{Z}}

\newtheorem{theorem}{Theorem}[section]
\newtheorem{lemma}{Lemma}[section]

\newtheorem{definition}{Definition}[section]

\newtheorem{remark}{Remark}[section]

\begin{document}
\title{\Large \bf Horizontal mean curvature flow as a scaling limit  of a mean field equation in the Heisenberg  group.
}

\author{Giovanna Citti \thanks{Università di Bologna, Bologna, Italy, e-mail:giovanna.citti@unibo.it} \and
 Nicolas Dirr\thanks{Cardiff School of Mathematics, Cardiff University, Cardiff, UK, e-mail: dirrnp@cardiff.ac.uk.} \and
Federica Dragoni
\thanks{Cardiff School of Mathematics, Cardiff University, Cardiff, UK, e-mail: DragoniF@cardiff.ac.uk.} \and
Raffaele Grande 
\thanks{The Czech Academy of Sciences, Institute of Information Theory and Automation, Prague, Czech
Republic, e-mail: grande@utia.cas.cz)}
}

\maketitle
\begin{abstract}
\noindent
We derive curvature flows in the Heisenberg group by formal asymptotic expansion of a nonlocal mean-field equation under the anisotropic rescaling of the Heisenberg group. This is motivated by the aim of connecting mechanisms at a microscopic (i.e. cellular) level to macroscopic models of image processing through a multi-scale approach. The nonlocal equation, which is very similar to the Ermentrout-Cowan equation used in neurobiology,  can be derived from an interacting particle model.  As sub-Riemannian geometries play an important role in the models of the visual cortex proposed by Petitot and Citti-Sarti, this paper provides a mathematical framework for a rigorous upscaling  of models for the visual cortex from the cell level via a mean field equation to curvature flows which are used in image processing. From a pure mathematical point of view, it provides a new approximation and regularization of Heisenberg mean curvature flow. {Using the local structure of the roto-translational group, we extend the result to cover the model by Citti and Sarti. Numerically, the parameters in our algorithm interpolate between solving an Ementrout-Cowan type of equation and a  Bence–Merriman–Osher algorithm type algorithm for sub-Riemannian mean curvature. We also reproduce some known exact solutions in the Heisenberg case.}
\end{abstract}

\begin{center}
\textbf{Keywords}: Mean Curvature flow, Carnot groups, Heisenberg group, asymptotic expansion, subellitpic PDEs, Sub-Riemannian geometries, H\"ormander condition.
\end{center}

\section{Introduction}

The derivation of interface evolution laws, an important example being the evolution by mean curvature, from microscopic models has attracted the attention of mathematicians for a long time and the literature is too vast to give a complete overview here, so we restrict ourselves to what is
necessary to understand the motivation for the  scaling limit  considered in this paper. For a model on a small scale, frequently an equation of reaction-diffusion type is chosen, see e.g \cite{al}, \cite{cah}, \cite{sh}. Here the fast reaction drives the system on most of the domain towards two distinct spatially constant equilibria. The evolution of the transition layer, separating the zones occupied by either equilibrium, follows a surface evolution law on a larger time-space scale,  for example motion by mean curvature. This is a fully deterministic model that neither accounts for the discrete microstructure of matter nor possible random effects. In order to take into account these effects, one could add some form of spatially regularised white noise, or one could use an interacting particle model on the small scale. The latter approach has the advantage of allowing scope for modelling the microscopic structure of the material under consideration and to specify noise on this smallest possible scale, called henceforth "atomistic." In the simple case under consideration here, this is a ferromagnetic Ising model, that can be described as follows:  the "cells" are placed on an integer lattice and either active $(+1)$ or inactive $(-1).$ The state of the cell at $X\in \mathbb{Z}^N$ ($N$ is the space dimension) at time $t$ is denoted by $\sigma(t,X).$
They switch from active to inactive and vice versa with a rate which depends on the 
(possible anisotropic) averaged input from their neighbours, preferring to be aligned with them. Unless in the special case of zero temperature (see \cite{lac2014}), rigorous derivation of macroscopic evolution laws from Ising models has been an open problem for decades. The situation is different for local mean field models.

The jump at $X$ depends  on the state of the cell at $X$ and a 
 sigmoidal-type function of the following
average over many cells:
\begin{equation}
\label{sigmoidal-type}
\gamma^N \sum_{Y\in \Z^N}J(\gamma X- \gamma Y)\sigma(t,Y).
\end{equation}
Here $J$ is a smooth nonnegative interaction kernel, depending on the two cell positions $X$ and $Y$ in $\Z^N$, with unit mass and  compact support. We expect the results to extend to kernels with exponential decay like the heat kernel, but this is beyond the scope of this paper. The small parameter $\gamma$ is the inverse of the interaction range. 
 In the limit $\gamma\to 0$, the processes $\sigma(t,X)_{X\in \mathbb{Z}^N}$ are approximately independent, which allow to derive mathematically rigorous limit theorems.  While $\gamma\to 0$ is a mathematical fiction, albeit one backed up by numerical experiments, \cite{Kat2001},  in the time-independent case there are more precise statements regarding  the relationship between the mathematical limit $\gamma \to 0$ and the case for finite but long-range interactions (Lebowitz-Penrose limit), see \cite{pre2}.
Let us denote the activity  averaged over a block of size $\gamma^{-\alpha},$ $0<\alpha<1,$ (smaller than the interaction range) by $m^\gamma$, i.e.
$$
{ m^\gamma(t,X,\omega):=\left(\gamma^{\alpha} \right)^N
\sum_{|Y-X|_\infty<\frac{1}{2}\gamma^{-\alpha}}
\!\!\sigma(t,Y,\omega)},
$$
where $|Y-X|_{\infty}=\max_{i=1,\dots, N}|Y_i-X_i|$.
Then, using the approximate independence of cells at different sites,  by a law of large numbers $m^\gamma$ becomes deterministic as $\gamma\to 0$ and solves a nonlocal evolution equation.
The mesoscopic limit that we are working with is given by $\lim_{\gamma\to 0 }m^\gamma(t,\gamma^{-1}x) =
m(t,x)$ in probability.  The limit function $m$ solves
\begin{equation}\label{tba}
\frac{\partial}{\partial t}m=-m+\tanh(\beta J* m +a),
\end{equation}
where $J$ is the kernel introduced in \eqref{sigmoidal-type}, $\beta>1 $ is the inverse temperature and $a$ may be a function or a constant.
This has been extended by Katsoulakis and Souganidis, \cite{kat}, to cover anisotropic kernels.
If  $a=0$ and the nonlocal equation is rescaled diffusively, i.e. time with the square of the spatial rescaling, motion by mean curvature is obtained as  limit evolution law in a similar way as for the Allen-Cahn equation. Let us consider the equation solved by $m(\varepsilon^{-2}t,\varepsilon^{-1}x)$, that is 
$$
\frac{\partial}{\partial t}m=-\varepsilon^{-2}(m+\tanh(\beta J^\varepsilon *  m )).
$$
Here $J^{\varepsilon}$ is the rescaled kernel, i.e. $J^{\varepsilon}(x)=\varepsilon^{_N}J(\varepsilon^{-1} x)$. The $\varepsilon^{-2}$ on the right hand side forces the solution to stay close to the two spatially constant stable equilibria $\pm m_\beta.$ The unscaled equation has a spatially non-constant solution connecting these equilibria,
called $\overline m$ with $\lim_{x\to \pm\infty}\overline m(x)=\pm m_\beta.$ We make  the ansatz $m(t,x)=\overline m (d(x,M_{t})),$ where $d$ is the signed distance (see Section 2 for a definition)  and $M_t$ a set such that its boundary $\partial M_t$ is a hypersurface evolving in time. An asymptotic expansion (see \cite{pres}) shows that
in order to be a solution at highest order (for $a=0$ and the kernel invariant under rotations) the interface has to evolve according to
$$
V=\theta \, \kappa,
$$
where $V$ is the normal velocity of the surface, $\theta$ is a parameter which depends on the kernel $J,$ and $\kappa$ is the mean curvature {(for anisotropic kernels and small but  non zero forcing  $a$ see \cite{kat}).} {

So far the underlying geometric structure was the usual Euclidean $\mathbb{R}^n$. In this paper, however, we consider  the analogous evolution in the 1-dimensional Heisenberg group, which is  a step 2 Carnot group, that can be identified with $\mathbb{R}^{3}$ with the standard topology, but with a different metric structure which is not equivalent to the Euclidean $\mathbb{R}^{3}$ at any scale. 
Roughly speaking, the Heisenberg group is a non-commutative Lie group, where a Riemannian metric is defined on a 2-dimensional distribution $\mathcal{H}$.

An admissible curve is any  absolutely continuous curve whose tangent vector belongs at every points to such distribution. 
It is well-known that, any two points can be connected by at least one admissible curve (Chow's Theorem, see e.g. \cite{cho}) and therefore we can define a distance between two points in $\mathbb{R}^3$ as the infimum over the lengths of all  admissible curves connecting these to points, the so-called Carnot-Carathéodory distance. This distance is 1-homogeneous with respect to an anisotropic scaling $x=(x_1,x_2,x_3)\mapsto\delta_\lambda(x)=(\lambda x_1, \lambda x_2,\lambda^2 x_2).$ 
This anisotropic scaling reflects the fact that not all directions are equal. Moving away from the span of the vector fields is possible via only interchanging between the vector fields (commutators).

This leads to formulate equation \eqref{tba} in the Heisenberg group, where now $J$ is the heat kernel of the Heisenberg group, or a general nonnegative, smooth kernel of the form $J(x_1^2+x_2^2,|x_3|)$ with compact support.
The rescaled kernel needs to be adapted to the anisotropic rescaling of the Heisenberg group, that is
 $J^{\varepsilon}(x)=\varepsilon^{-Q}J\left(\delta_{\frac{1}{\varepsilon}}(x)
\right)$, where $Q$ is the homogeneous dimension, i.e. $Q=4$.
 The asymptotics of this heat kernel for small times are similar to the Euclidean heat kernel where the Carnot-Carathéodory distance replaces the Euclidean one.
With the mentioned anisotropic spatial scaling we obtain by a formal asymptotic expansion a similar evolution law, motion by Heisenberg mean curvature. Of course this relies on symmetry assumptions on the kernel which should be symmetric if restricted to the ``allowed'' submanifold.

In spite of the formal similarity, the Heisenberg mean curvature flow is very different from its Euclidean counterpart due to the existence of {\em characteristic points} where the interface  evolution is not defined.  Characteristic points, which are present in every compact surface, make this evolution far more difficult to understand than the standard Euclidean one; in fact to our knowledge only partial results for the comparison principles for the Heisenberg mean curvature flow are known, see \cite{{ba},{cap},{fer}}. This is why we do not strife for a rigorous convergence result as e.g. in \cite{kat}. The presently available techniques  do not allow to prove such a general theorem in the framework of viscosity solutions due to lack of  comparison principles. 

On the other hand, the nonlocal convolution equation for small but finite $\varepsilon$ is always defined, and the zero level set of solutions stays close to a surface evolving by Heisenberg mean curvature flow in the absence of characteristic points or other singularities.
This means that this paper presents a regularisation of Heisenberg mean curvature flow. Other approximations and regularisations have been proposed previously e.g. in  \cite{ba}, \cite{dirr}, but this is the first one rooted in a multiscale analysis which can in principle be related to cellular models.

Numerical simulations show that the zero level set of solutions to our nonlocal evolution equation is a good approximation to evolution by Heisenberg mean curvature flow even in the presence of characteristic points, see 
Section  4.2 for a situation where a unique exact solution  of the mean curvature flow with characteristic points is known.

Later we exploit the local structure of the roto-translational geometry SE(2) (see Section 5 for a definition) in order to extend the results from Heisenberg to SE(2). The Heisenberg group is the tangent cone to SE(2), in the sense that its Lie algebra provides an approximation of the Lie algebra of SE(2) up to step 2 commutators; then by the exponential map we can obtain the corresponding local approximation of the two groups (see Rothschild-Stein freezing procedure \cite{RS}).

Our approach is motivated by the models of the visual cortex proposed by J. Petitot and Y. Tondut, \cite{PT},  and by the Citti-Sarti model for the visual cortex, \cite{cisar2}, \cite{pe}. Cells are activated by receptors corresponding to light reaching a (discrete) point $(x,y)$ on the retina. The visual signal is carried to the cortex by the optical nerve, and thanks to the cortical connectivity, simple cells of the primary cortex V1 are able  to extract information on the direction of the gradient of the light intensity, which can be expressed by an angle $\theta.$ 

So the two-dimensional image is ``lifted'' to a three-dimensional cortical representation, but with constraints due to the fact that it was a lift of a two dimensional image. In \cite{PT}, the model was expressed as a constraint minimization problem in the Heisenberg group, where the problem can be locally described and in \cite{cisar} as a sub-Riemannian structure on $\mathbb{R}^2 \times S^1$, the roto-translational geometry, which allows a global description of the problem. The Heisenberg group is the tangent cone to $\mathbb{R}^2 \times S^1$ so that the two geometries are locally equivalent.  Again, two vector fields coming from a Lie-group structure span all admissible directions and Chow's theorem holds, a Carnot-Carath\'eodory distance can be defined. Locally the structure is the one of the Heisenberg group.

In \cite{cisar} two neural mechanisms of diffusion along the cortical connectivity and non maximal suppression where modelled as a sub-Riemannian mean curvature flow able to perform for image completion and denoising the corresponding sub-Riemannian mean curvature flow can be used for image completion and denoising. A mesoscopic description of cell activity going back to Ermentrout and Cowan \cite{co}, yields a nonlocal evolution equation like \eqref{tba}, and in   Bressloff and Cowan, \cite{BC},  a kernel $K$ invariant with respect to the geometry of SE(2) was proposed. Our approach allows in principle  to connect a cellular description first  to a mesoscopic Ermentrout-Cowan description and then all the way up the scales to mean curvature flow,  providing a rigorous link between descriptions by nonlocal equations and descriptions by geometric evolution equations.
 \\

The paper is organised as follows: in the Section 2 we recall briefly some results about the mean field equation and its connection with the mean curvature flow in the Euclidean setting, and we define our interaction kernel in the Heisenberg group. In the Section 3 we will give some definitions related to sub-Riemannian geometries and the horizontal mean curvature flow. In the Section  4.1 we  study the asymptotic expansion at non characteristic points and in Section 4.2 we report results on numerical experiments, in particular consistency for an exact solution from \cite{fer} which has characteristic points. In Section 5 we extend the theoretical results from Section 4 to the roto-translational geometry.

\section{The Euclidean case}
We recall in slightly more detail how one obtains motion  by mean curvature from the nonlocal mean field equation \eqref{tba}:
 it is in the Euclidean setting given by  
\begin{equation} \label{mer}
\frac{\partial}{\partial t}m + m - \tanh \big[ \beta \; ( J * m +a)\big] = 0, \ \ \ \mbox{in} \ \ [0, \infty ) \times \mathbb{R}^{N}  ,
\end{equation}
where $J: \mathbb{R}^{N} \rightarrow [0,\infty)$ is a function called \emph{kernel} with compact support and it has form $J(x)=J( | x | ),$ $\int_{\mathbb{R}^N}J(x)dx=1$  and $\beta>1$, $a \in \mathbb{R}$. Note that we have generalised the setting slightly by introducing a constant force $a$.

It is important to note that, since $J \geq 0$, the Equation \eqref{tba} admits a comparison principle between solutions. 
\begin{lemma} (\cite{kat}, \cite{pres})\label{comparison}
 If $w_{1}$ and $w_{2}$ solve \eqref{mer}, then
\begin{equation*} 
 w_{1} \leq w_{2} \ \ \mbox{on} \ \ \{ 0 \}  \times \mathbb{R}^{N}  \ \Rightarrow \  \ w_{1} \leq w_{2} \ \mbox{on} \  (0, \infty)  \times \mathbb{R}^{N}.
\end{equation*}
\end{lemma}
This holds for any non-negative kernel $J(x,y)$, regardless of symmetries, so it carries directly over to our sub-Riemannian setting.
\begin{remark}
There exists $a_0>0$ such that the mean field equation \eqref{mer}  allows for $|a|<a_0$ three spatially constant steady states $m^{-}_{\beta , a } <  m^{0}_{\beta, a } < m^{+}_{\beta , a }$ (\cite{pres}) i.e. constant solutions of the algebraic equation
\begin{equation} \label{eq1}
m = \tanh\left[\beta \left(\left  (\int_{\mathbb{R}^N}J(x)dx\right) m +  a\right)\right].
\end{equation}
This follows from
\begin{equation*}
\beta \int_{\mathbb{R}^N}J(x)dx > 1.
\end{equation*}
In the case $ a = 0$
\begin{equation*}
m^{\pm}_{\beta, 0} = \pm m_{\beta} \ \ \mbox{and} \ \ m^{0}_{\beta, 0} = 0.
\end{equation*} 
Furthermore, it is immediate to remark that, by using $\int_{\mathbb{R}^N}J(x)dx=1$, Equation \eqref{eq1} may be rewritten as
\begin{equation*}
\frac{1}{\beta} {\rm arctanh} (m) = m+a.
\end{equation*}
\end{remark}
Now let us consider the real line, i.e. $N=1$, later we will  reduce a general kernel to a one-dimensional one by integrating out the additional degrees of freedom in a suitable way. Let us denote the 1-dimensional kernel by $\overline J(r),$ $r \in \mathbb{R}$.
Consider Equation \eqref{mer} replacing the forcing  $a$ by $\varepsilon a$ for $\varepsilon$ sufficiently small. In this case we have  so called \emph{travelling waves}, which are solutions of \eqref{mer} in one dimension  and have the form
\begin{equation*}
m^\varepsilon(t,r) = q(r-c^{\varepsilon}_at;{\varepsilon a}).
\end{equation*}
$q$ solves 
\begin{equation*}
\left\{
\begin{aligned}
&\quad c^{\varepsilon}_a \dot{q}( r ;\varepsilon a) + q( r; \varepsilon a) = \tanh \beta [\overline{J} * q(r; \varepsilon a) + \varepsilon a ],\\
&\quad  \dot{q} >0,\\
& \lim_{r\to \pm \infty}q(r, \varepsilon a) = m^{\pm}_{\beta , \varepsilon a}.
\end{aligned}
\right.
\end{equation*}

Note that the limit has  exponential rate of convergence, see \cite{pres} for details.

The case $a=0$ is particularly important. In this case we have a stationary solution which is not constant in space, called, following \cite{pres, pre2}, an {\em instanton} and denoted henceforth by $\overline m.$  For further reference, we summarise its properties in  the following lemma.

\begin{lemma}(\cite{pres, pre2})\label{instantonlemma} Suppose the kernel $\overline J:\mathbb{R} \to \mathbb{R}$ with  $\overline  J(r)= \overline  J(|r|)$ is smooth, compactly supported, nonnegative
and $\|\overline  J\|_{L^1}=1,$ and moreover $\beta>1$.  
Then there exists a unique  function $\overline m: \mathbb{R} \to (-m_\beta,+m_\beta)$ such that 
\begin{equation}\label{insteq}
-\overline m+\tanh\left(\beta \overline J\ast \overline m\right)=0\quad {\rm in\ } \mathbb{R},
\end{equation} and
\begin{itemize}
\item $\overline m(0)=0;$
\item $\overline m(-r)=-\overline m(r);$
\item $\overline m$ is smooth;
\item $\overline m^\prime(r)=\frac{d}{dr}\overline m(r)>0$;
\item $\lim_{r\to\pm \infty}\overline m(r)=\pm m_\beta$ with exponential rate of convergence;
\item  $\lim_{r\to\pm \infty}\overline m^{(k)}(r)=0$ with exponential rate of convergence and where $\overline m^{(k)}$ is the derivatives of order $k\ge 1$.
\end{itemize}
\end{lemma}
 For constructing an asymptotic expansion both in the Euclidean and sub-Riemannian case, the linearisation of Equation \eqref{insteq} around $\overline m$ is important.
\begin{lemma}(\cite{pres})\label{linop}
Given $\overline m$ as in Lemma \ref{instantonlemma}, we define the linear operator on continuous functions $f: \mathbb{R} \to \mathbb{R}$ with compact support as
\begin{equation}
\label{L}
(\mathcal L f )(r)=\frac{\beta}{1-(\overline m(r))^2}(\overline J \ast f)(r) -f (r),
\end{equation}
then it extends to a  self-adjoint operator on $L^2(\mu)$ with $\mu(dx)=\frac{1}{1-(\overline m(r))^2}dx.$ This operator has nonnegative real spectrum, an eigenvalue zero with eigenfunction $\overline m ^\prime$ and a spectral gap, i.e. the nonzero part of the spectrum is contained in some interval $(\tau,\infty)$ with $\tau>0.$ As a consequence for all $g\in L^2(\mu)$, the equation
$$
(\mathcal L f )(r)=g(r)
$$ has a solution in $L^2(\mu)$ if and only if 
$$
0=\langle g, \overline m ^\prime\rangle_{L^2(\mu)}=\int_\mathbb{R} \frac{\overline m^\prime (r)}{1-(\overline m(r))^2}g(r)dr.
$$
\end{lemma}

This function  $\overline m $ describes the transition from the $-$ to the $+$ phase. 
Our aim is to convergence under diffusive rescaling to  motion by mean curvature, so let us  recall briefly 
the \emph{generalized evolution} of an hypersurface, with velocity $V= \theta \kappa + A$ where $\theta> 0$, $A \in \mathbb{R}$ and $\kappa$ is the mean curvature.

The Euclidean signed distance  $d_E(x,S)$ of a point  $x\in \mathbb{R}^N$ from a set $S\subseteq \mathbb{R}^N$ is defined as
$$
d_E(x,S)=\begin{cases} -d(x, \mathbb{R}^N\setminus S),  \quad \quad \mbox{if} \  x \in S, \\ d(x, S),  \quad \quad \quad \quad \ \ \ \ \mbox{if} \ x\in \mathbb{R}^N\setminus S. \end{cases} \label{dista}
$$ 

\begin{definition}
Let $M_{0} \subset \mathbb{R}^{N}$ be open and bounded and we denote by
$
\partial M_{0}
$ the boundary of $M_0$. (Note that $x\in \partial M_0 \; \Longleftrightarrow\; d(x,\R^N\backslash M_0)= d(x,M_0)=0$.)
The \emph{generalised evolution} of the  hypersurface $\partial M_0$ is defined as the hypersurface  $\partial M_{t}$  where 
\begin{equation*}
M_{t} := \{ x \in \mathbb{R}^{N}  : \ w(t,x)>0 \},
\end{equation*}
and $w$ is the unique viscosity solution of
\begin{equation} \label{zacc}
\begin{cases} w_{t} - \theta\, tr \left( I - \frac{Dw \otimes Dw}{|Dw|^{2}} \right) D^{2} w - A |Dw| = 0, \ \ \mbox{in} \ (0, \infty)  \times \mathbb{R}^{N}, \\ w= d_{0}, \ \ \  \ \ \  \ \ \  \ \ \  \ \ \  \ \ \  \ \ \  \ \ \  \ \ \  \ \ \  \ \ \  \ \ \  \ \ \  \  \ \ \ \ \ \mbox{on}   \ \{ 0 \}  \times \mathbb{R}^{N}, \end{cases}
\end{equation}
and 
\begin{equation}\label{d0}
d_{0}(x) = d_E(x,M_0).
\end{equation}
\end{definition}
Roughly speaking, under diffusive rescaling with $\varepsilon$ small, solutions will be close to $\pm m_\beta$ except near a moving surface $\partial M_t$, where they will look like $\overline m(\varepsilon^{-1}d_E(x,M_t)),$ where $d_E(x,M_t)$ is the signed distance.
A precise definition of the rescaling is as follows:
\begin{remark}
If $m$ solves \eqref{mer}, then 
\begin{equation*}
m^{\varepsilon}(t,x)= m(  \varepsilon^{-2} t,  \varepsilon^{-1} x), \ \ (t,x) \in  (0, \infty) \times \mathbb{R}^{N} 
,
\end{equation*}
solves 
\begin{equation} \label{tir}
m_{t}^{\varepsilon} + \varepsilon^{-2}[m^{\varepsilon} - \tanh \beta (J^{\varepsilon} * m^{\varepsilon} +  a)] =0, \ \ \mbox{in} \ \ (0, \infty)  \times \mathbb{R}^{N},
\end{equation} 
with $
J^{\varepsilon}(x) = \varepsilon^{-N} J( \varepsilon^{-1} x)$, and $x \in \mathbb{R}^{N}.$
This rescaling is immediate since the Euclidean space is \emph{isotropic}, it has to be replaced by suitable dilations  in the case of a sub-Riemannian geometry (see Section 3 for further details).
\end{remark}
It is possible to find a coefficient  $\theta$ (see \cite{pres, kat}) such that solutions to the scaled Equation \eqref{tir} converge
to the motion by mean curvature in the following sense: let us suppose that $M_t$ is a set such that its  boundary $\partial M_t$ which evolves by the law 
\begin{equation} \label{talk}
V= \theta \kappa.
\end{equation}
 In this case $M_t$ may be the set $\{w>0\}$ with $w$  as in \eqref{zacc}, or we may have a classical evolution by mean curvature flow of the smooth surface $\partial M_t.$ 
 Furthermore we assume that the kernel $J$ is rotationally symmetric, i.e. depends on $x\in \mathbb{R}^N$ only through its Euclidean norm $|x|$, thus we introduce   the one dimensional kernel 
 $
\overline J:\R\to [0,+\infty)$ as
$$
\overline J (r)=\int_{\mathbb{R}^{N-1}}J(r, x_2,\ldots,x_N)dx_2\dots dx_N.
$$ 
Note that $\overline J$ satisfies the assumptions of Lemma \ref{instantonlemma}, then the following result hold true.
\begin{theorem}[\cite{pres, kat}] \label{bod}
Let $m^{\varepsilon}$ be the solution of (\ref{tir}) with $a=0$ and with initial datum 
$$
m^\varepsilon(0,x)=\overline m\left(\frac{d_0(x)}{\varepsilon}\right)
$$ 
where $d_0$ is as in \eqref{d0} and $\overline m$ is as in Lemma \ref{instantonlemma}. Then, as $\varepsilon \rightarrow 0^{+}$
\begin{equation*}
m^{\varepsilon} \rightarrow \begin{cases} m^{-}_{\beta, 0}, \ \ \ \ \ \ \ \ \ \ \  {\rm in\ } \mathring  M_t, \\ m_{\beta, 0}^{+}, \ \ \ \ \ \ \ \ \ \ \  {\rm in\ } \mathbb{R}^N\setminus \overline{M_t}.\end{cases}
\end{equation*}
locally uniformly in $\R^n\backslash \partial M_t$.
\end{theorem}

The main idea is the ansatz for a solution
\begin{equation}\label{ansatzeucl}
m^\varepsilon(t,x)=\overline m \left(\frac{d_E(x,  M_t)}{\varepsilon}\right)+\varepsilon m_1 \left(\frac{d_E(x,  M_t)}{\varepsilon}\right),
\end{equation}where $m_1$ is to be determined, and we suppose to have a classical solution for the mean curvature flow Equation \eqref{zacc}.

Substituting ansatz \eqref{ansatzeucl} in the Equation \eqref{tir} and writing $r=\varepsilon^{-1}d_E(x,  M_t)$, by a second order Taylor expansion of the distance function and a first order Taylor expansion of the hyperbolic tangent around $r$, we get
\begin{equation}
\label{firstequation}
\begin{aligned}
\varepsilon^{-1}\frac{\partial}{\partial  t} d (x,  M_{t})\overline m^\prime(r)&=
\varepsilon^{-2}\left[\tanh\left(\beta (J\ast \overline m)(r)\right)-\overline m(r)\right]\\
&+\varepsilon^{-1}\left[\mathcal {R}(t,r,x)+\mathcal L m_1(r)\right]+O(1),
\end{aligned}
\end{equation}
where the linear operator $\mathcal L$ is given in \eqref{L} and $\mathcal R$ is defined by
the expansion of  the nonlinear term in \eqref{tir} as $\tanh\left(\beta (J\ast \overline m)(r)\right)-\overline m(r)+\varepsilon
{\mathcal R}(t,r,x)+O(\varepsilon^2)$ 
First note that the term of order $\varepsilon^{-2}$ on the right hand side of \eqref{firstequation} vanishes because of \eqref{insteq}.
If $d_E(x,  M_t)\gg \varepsilon$ {\em all} terms are exponentially small due to the exponential convergence of the instanton towards its limits at $r\to \pm \infty,$ (see Lemma \ref{instantonlemma}) so the equation is solved up to an exponentially small error. Near the interface, 
at highest order, $\frac{\partial}{\partial  t} d (x,  M_{t})=V,$ the normal velocity, and the term $\mathcal R$ can be expressed as $\mathcal 
\kappa(\tilde x)R(r),$ where $\kappa$ is the mean curvature and $\tilde x$ the point on $\partial M_t$ closest to $x$, plus  highest order terms.
Now we can split $R(r)$ in its projection on $\overline m$ with respect to the $L^2(\mu)$-norm and in a part which is orthogonal to 
$\overline m^\prime$, i.e.
$$
R(r)=\theta \overline m ^\prime + R^\perp(r),\quad \langle R^\perp(r),\overline m^\prime\rangle_{L^2(\mu)}=0,
$$where this equation defines the coefficient $\theta$ in the evolution by mean curvature. 
If we insert this splitting in
\eqref{firstequation}, we see that the part multiplying $\overline m^\prime$ cancels with the part on the left hand. As we can find $m_1$ such that 
$
\mathcal L m_1(r)+R^\perp(r)=0
$ by Lemma \ref{linop}, the Equation \eqref{tir} is solved by our ansatz up to order 1 terms  near the interface and up to exponentially  small terms far from the interface. For details we refer to \cite{pres}.
If the limit surface evolution law has a classical solution, then in principle it is possible to find functions 
$m_k(r),\ k\ge 1$ and corrections to the distance function $d_k$ 
such that 
$$
m^\varepsilon(t,x)\!=\!\overline m \left(\frac{d_E(x,  M_t)+\sum_{k=1}^K \varepsilon^k d_k(x)}{\varepsilon}\right)\!+\!
\sum_{k=1}^K\varepsilon^k m_1 \left(\frac{d_E(x,  M_t)+\sum_{k=1}^K \varepsilon^k d_k(x)}{\varepsilon}\right)
$$
solves the equation up to an error of sufficiently small order to bound the difference between this ansatz and the actual solution in the $L^2$ norm. For the Allen-Cahn equation this was done in \cite{sh}, and for nonlocal equations similar to \eqref{tir} e.g. in \cite{enza, dirr2}.
In the presence of a comparison principle both for \eqref{tir} and the limit evolution, it is not necessary to construct any further terms in the expansion. Instead we insert \eqref{ansatzeucl} in \eqref{tir} biased by $a=\pm\delta\varepsilon$ for some $\delta>0.$ We can choose $\delta=\delta(\varepsilon)$ with $\lim_{\varepsilon\to 0} \delta(\varepsilon)=0$ such that the  the ansatz becomes a subsolution/supersolution for the Equations \eqref{tir} with 
$a=\pm\delta\varepsilon.$ Hence both our ansatz \eqref{ansatzeucl} and the exact solution of \eqref{tir} with $a=0$ stay between the solutions for $a=-\delta \varepsilon$ and $a=+\delta \varepsilon.$ It remains  to show that the zero level sets of these two biassed solutions stay close. By a formal argument as the none following \eqref{firstequation} we can find an effective evolution law for the limit surface of the form $V=\theta \kappa\pm c(\delta),$ with $c$ monotone and $c(0)=0.$ In this way both the zero level set of the actual solution of \eqref{tir} with $a=0$  and of our ansatz is trapped between two interfaces that stay close to each other and close to $\partial M_t$ as $\varepsilon \to 0.$  

This  reasoning, based on the use of comparison principles, can be extended to a viscosity proof, see \cite{kat}. 
In the case of an anisotropic kernel (see \cite{kat}) the instanton and the travelling waves depend in addition on a direction, the direction of the surface normal.

We would like to extend this ansatz to a different geometry, that of the Heisenberg group, and to its anisotropic scaling. Let us recall some basic notions.

\section{The Heisenberg group}
In this section we  recall some properties related to the 1-dimensional Heisenberg group $\mathbb{H}^1$.  This is the simplest but most important example of a sub-Riemannian geometry. For a general definition of sub-Riemannian geometries and their properties we refer to \cite{mont}. The Heisenberg group is a step 2 Carnot group, i.e.  a simple connected, nilpotent  non-commutative Lie group: details on the Heisenberg group and Carnot groups can be found e.g. in \cite{BLU}. 
Let us start by recall that the  Baker–Campbell–Hausdorff formula in this case writes simple as
 \begin{equation*}
 \exp(X_{1}) \exp(X_{2}) = \exp \bigg(X_{1} + X_{2} + \frac{1}{2}[X_{1}, X_{2}] \bigg),
 \end{equation*}
where the vector fields can be chosen as
\begin{equation*}
X_1(x)= \begin{pmatrix}
1\\ 0 \\ - \frac{x_{2}}{2}
\end{pmatrix}
\quad  \textrm{and} \quad
X_2(x)=
\begin{pmatrix}
0\\1 \\ \frac{x_{1}}{2}
\end{pmatrix},
\quad \forall \; x=(x_1,x_2,x_3)\in \mathbb{R}^3.
\end{equation*} 
$X_1$ and $X_2$ are the left-invariant vector fields w.r.t. to the following non-commutative group law:
\begin{equation*}
(x_{1}, x_{2}, x_{3}) \circ (y_{1}, y_{2}, y_{3}) = \bigg(x_{1}+y_{1}, x_{2} + y_{2}, x_{3}+y_{3} + \frac{1}{2}\big(x_{1}y_{2} - x_{2}y_{1}\big) \bigg).
\end{equation*}
We will sometimes denote the left-translations as $L_{x}(y):= x \circ y$.

Note also that, for any $x\in \mathbb{R}^3$, $$X_{3} (x)=[X_1,X_2](x)=\begin{pmatrix}0\\0\\1\end{pmatrix},$$ hence the H\"ormander condition is satisfied. By applying 
 the  \emph{Chow's Theorem} (see e.g. \cite{cho}), it is always possible to connect any two given points  by a horizontal curve, i.e. by an absolutely continuous curve $\gamma:[0,T] \rightarrow \mathbb{R}^{3}$  such that 
\begin{equation*}
\dot{\gamma}(t) =  \alpha_{1}(t) X_{1}(\gamma(t))+ \alpha_{2}(t) X_{2}(\gamma(t)), \ \ \mbox{for a.e} \ t \in [0,T],
\end{equation*}
where $\alpha_{1}$ and $\alpha_{2}$ are two scalar measurable functions.
The Carnot-Carathéodory distance, between two points $x$ and $y$, 
 $d: \mathbb{H}^{1} \times \mathbb{H}^{1} \rightarrow \mathbb{R}$ is defined by
\begin{equation*}
d(x,y)= \inf\{l(\gamma)| \mbox{$\gamma$ horizontal curve connecting $x$ and $y$}\},
\end{equation*}
where 
$l(\gamma) := \int_{0}^{T}  \sqrt{\alpha_{1}^{2}(t) + \alpha_{2}^{2}(t)} dt$. \\
Given $M \subset \mathbb{H}^{1}$ open, with boundary $\partial M$, we define the distance from the boundary  as $d(x,  M) = \inf_{y \in \partial M} d(x,y)$. The \emph{signed distance}  is 
\begin{equation} \label{sig}
d_{\mathbb{H}^{1}}(x,  M)= \begin{cases} d(x,  M), \quad \ \ \ x \in \mathbb{H}^{1} \setminus M, \\ - d(x, M),  \quad \ x \in M. \end{cases}
\end{equation}

Let us recall  that $\mathbb{H}^{1}$ is \emph{anisotropic},  in fact the dilatations are defined  as
\begin{equation*}
\delta_{\lambda}(x_{1},x_{2},x_{3})= (\lambda x_{1} , \lambda x_{2} , \lambda^{2} x_{3}), \ \ \mbox{for $\lambda\geq 0$.}
\end{equation*}
The vector fields associated to the Heisenberg group allows us to define the horizontal version of some classical differential operators. Let $u: \mathbb{H}^{1} \rightarrow \mathbb{R}$ be a $C^2$ function, we can define the horizontal gradient as
\begin{equation*}
\nabla_{\mathbb{H}^{1}} u =(X_{1} u ,  X_{2} u)^{T}, 
\end{equation*}
and the symmetrized horizontal Hessian as
\begin{align*}
(D^{2}_{\mathbb{H}^{1}} u)^{*} = \begin{bmatrix} X_{1}(X_{1}u) & \frac{X_{1}(X_{2}u) + X_{2}(X_{1}u)}{2} \\ \frac{X_{1}(X_{2}u) + X_{2}(X_{1}u)}{2} & X_{2}(X_{2}u)   \end{bmatrix}.
\end{align*}
 Here the horizontal Laplacian is
\begin{equation*}
\Delta_{\mathbb{H}^{1}}u = X_{1}(X_{1}u) + X_{2}(X_{2}u),
\end{equation*} 
For $v: \mathbb{H}^{1} \rightarrow \mathbb{R}^{2}$ $C^1$-function,
 the horizontal divergence is given by
\begin{equation*}
div_{\mathbb{H}^{1}} v(x) = (X_{1} v_{1}) + (X_{2} v_{2}).
\end{equation*}
Given $M$ as above, the \emph{horizontal normal}  $n_{0}(x)$ at the point $x \in \partial M$ is the renormalized projection of the Euclidean normal $n_E(x)$ on the horizontal space $\mathcal{H}_x=\textrm{Span}\big(X_1(x), X_2(x)\big)$.
The  \emph{horizontal mean curvature} is defined as the horizontal divergence of the horizontal normal, i.e.
$
k_{0}(x) = div_{\mathbb{H}^{1}} n_{0}(x).
$

If the projection of the Euclidean normal at the point $x \in \partial M$ vanishes, then $x$ is called a \emph{characteristic point}. Examples of characteristic points in $\mathbb{H}^{1}$ are the north and south poles of a sphere  centred at the origin, see e.g. \cite{dirr}.
\begin{definition} \label{c6glo}
Let $ \{ \partial M_t \}_{t \geq 0}$ be a family of smooth hypersurfaces in $\mathbb{H}^1$, depending on a time parameter $t\geq 0$. 
We say that $ \{ \partial M_{t} \}_{t \geq 0}$ is an \emph{evolution by horizontal mean curvature flow} 
of some hypersurface $\partial M$ if and only if $\partial M_0=\partial M$ and for any smooth horizontal  curve $\gamma: [0,T] \rightarrow \mathbb{H}^{1}$ such that $\gamma(t) \in \partial M_{t}$ for all $t \in [0,T]$, 
the horizontal normal velocity $v_{0}$ satisfies
\begin{equation} \label{sla}
v_{0}(\gamma(t)) = - k_{0}(\gamma(t)) n_0(\gamma(t)), \ \ t \in [0,T],
\end{equation}
where $n_{0}$ and $k_{0}$ as respectively the horizontal normal and the horizontal mean curvature.
\end{definition}
Equation \eqref{sla} is not defined at characteristic points, therefore to study the evolution there is necessary to use generalised approaches as  to study singularities in the standard Euclidean evolution, e.g. the so called \emph{level set approach} that was developed by Chen, Giga and Goto \cite{che} and Evans  and Spruck \cite{ev} for the Euclidean evolution and generalised to  the sub-Riemannian setting in
\cite{{ba},{cisar}, {cisar2}, {dirr},{fer},{gra3}}. 

Let us recall a theorem which connects the signed distance of $\mathbb{H}^{1}$ with the horizontal mean curvature of an hypersurface.  
\begin{theorem}(\cite{arc}) \label{eq}
Let $ M \subset \mathbb{H}^{1}$ be open and such that $\partial M \in C^{3}$, at any non characteristic point $x \in \partial M$, 
\begin{equation*}
\Delta_{\mathbb{H}^{1}} d_{\mathbb{H}^{1}}(x, M) = k_{0}(x),
\end{equation*}
where $d_{\mathbb{H}^{1}}$ is as in \eqref{sig} and $k_{0}$ the horizontal mean curvature of $\partial M$.
\end{theorem}
Furthermore, in the same paper, in the case of  graphs, the authors give an explicit formula at non characteristic points for the horizontal mean curvature flow and a property about characteristic points which will be useful later.
\begin{theorem}(\cite{arc}) \label{chp}
Let $M \subset \mathbb{H}^{1}$  be open and locally defined by the inequality $z - f(x,y)<0$, where $f$ is a $C^2$-function, and introduce
\begin{equation*}
E=  X_1X_1f (X_2f)^2 - 2(X_1X_2)^{*}f X_1f X_2f + X_2X_2f (X_1f)^2,
\end{equation*}
where $(X_1X_2)^{*}f= \frac{1}{2} (X_1X_2+ X_2X_1)f$. \\
Then, the symmetrized horizontal Hessian of the signed distance $d_{\mathbb{H}^{1}}$, at any non-characteristic point, is given by
\begin{equation*}
(D^2_{\mathbb{H}^{1}} d_{\mathbb{H}^{1}})^{*}= \begin{bmatrix} \frac{(X_2f)^2}{|\nabla_{\mathbb{H}^1} f|^5} E + 4\frac{X_1f X_2f}{|\nabla_{\mathbb{H}^1} f|^3} & -\frac{X_1f X_2f}{|\nabla_{\mathbb{H}^{1}} f|^5} E + 2 \frac{(X_1f)^2 - (X_2f)^2}{|\nabla_{\mathbb{H}^{1}} f|^3} \\  -\frac{X_1f X_2f}{|\nabla_{\mathbb{H}^{1}} f|^5} E + 2 \frac{(X_1f)^2 - (X_2f)^2}{|\nabla_{\mathbb{H}^{1}} f|^3}  & \frac{(X_1f)^2}{|\nabla_{\mathbb{H}^1} f|^5} E - 4 \frac{X_1f X_2f}{|\nabla_{\mathbb{H}^1} f|^3} \end{bmatrix},
\end{equation*}
and the horizontal Laplacian by
\begin{equation*}
\Delta_{\mathbb{H}^{1}}d_{\mathbb{H}^{1}} = \frac{E}{|\nabla_{\mathbb{H}^1}f |^{3}} = \frac{X_1X_1f(X_2f )^{2} - 2(X_1X_2)^{*}f X_1f X_2f + X_2X_2f (X_1f )^{2}} {| \nabla_{\mathbb{H}^1}f|^3}.
\end{equation*}
If $x_0$ is a characteristic point in $\partial M$, then $|D^2_{\mathbb{H}^1} d_{\mathbb{H}^{1}}(x,  M)|\! = \! O(d_{\mathbb{H}^{1}}(x,x_0)^{-1})$ as \\ $x \rightarrow x_0$. 
\end{theorem}

\section{Main theorem}
Let us make the following assumptions on the kernel.
\begin{definition}\label{kernel}
Let $J:\mathbb{R}^3\to[0,\infty)$ be such that 
\begin{itemize}
\item $J$ is smooth and compactly supported;
\item $J$ is radially symmetric in the first two variables and symmetric in the third variable, i.e. 
$J(x_1,x_2,x_3)=J(x_1^2+x_2^2,|x_3|).$
\end{itemize}
We define the following dimensionally reduced kernels,
\begin{eqnarray*}
\hat J (x_1,x_2)&=&\int_\mathbb{R} J(x_1,x_2,x_3)dx_3,\\
\overline  J (r)&=&\int_{\mathbb{R}^2} J(r,x_2,x_3)dx_2 dx_3,
\end{eqnarray*}
and the rescaled kernel as
$$
J^\varepsilon(x_1,x_2,x_3)=\varepsilon^{-4}J(\delta_{\varepsilon}(x_1,x_2,x_3))=\varepsilon^{-4}J(\varepsilon^{-1}x_1,\varepsilon^{-1}x_2,\varepsilon^{-2}x_3)
$$
Note that the power of the rescaling is given by  the homogeneous dimension that in $\mathbb{H}^1$ is equal to 4.
Let $\overline m$ denote the instanton with kernel $\overline J$ as in  Lemma \ref{instantonlemma}.
\end{definition}
Note that our definition of $\overline J$ relies on the radial symmetry of the kernel. In the general case, one would have to 
integrate over a direction in the $(x_1,x_2)$-plane and the result would depend on this direction, but by radial symmetry we may choose the second unit vector as direction.

Now we can define the nonlocal evolution equation in the Heisenberg group. We assume $m^\varepsilon$ solves
on $[0,T] \times \mathbb{R}^3$ the following:
\begin{equation}\label{Hnonloc}
\left\{
\begin{aligned}
&\varepsilon^2\frac{\partial}{\partial t}m^\varepsilon(t,x)= -m^\varepsilon(t,x)+\tanh\left(\beta\int_{\mathbb{R}^3} J^\varepsilon(y^{-1}\circ x) m^\varepsilon(t,y) dy\right),\\
& m^\varepsilon(0,x)=\overline m\left(\frac{d_{\mathbb{H}^1} (x,M_0)}{\varepsilon}\right).
\end{aligned}
\right.
\end{equation}
Here $M_0$ is supposed to have a $C^3$-boundary with no characteristic points.
Let us state our assumptions on the surface evolution.

\begin{definition}\label{HMCFDef}
Suppose  that, for some $T>0$, and some $\delta_0>0$, there exists for any $\delta\in (-\delta_0,\delta_0)$ a family of sets $M_t^\delta$ such that $M_0^\delta$ has Hausdorff distance $\delta $ from $M_0,$ $M_t^\delta $ has a smooth boundary without characteristic points for $t\in [0,T]$ and $\partial M_t^\delta$ satisfies $$V_0=\theta\kappa_0+
\delta$$
 in the classical sense with $\theta>0$ that will be chosen later   in \eqref{theta}. Moreover, if $\delta_1\le \delta_2$ and $M_0^{\delta_1} \subseteq M_0^{\delta_2},$ then  $M_t^{\delta_1} \subseteq M_t^{\delta_2}$ on $[0,T].$ (Comparison principle).  \\
In addition, the Hausdorff distance of $\partial M_t^\delta$ from $\partial M_0^\delta$ tends to zero as $\delta\to 0$, uniformly in 
$t.$
\end{definition}
Note that while  in the Euclidean case this is always true for small times provided $\partial M_0$ is smooth, in the Heisenberg case this is an assumption, since characteristic points may act as singularities, even if the surface is smooth. 

Now we are able to state our main theorem.

\begin{theorem}\label{mainthm}
Suppose the assumptions in Definition \ref{HMCFDef} are satisfied, $J$ satisfies Definition \ref{kernel}  for $\delta=0$ and $\theta$ as in \eqref{theta} and $m^\varepsilon$ solves
\eqref{Hnonloc}, then 
\begin{equation*}
m^{\varepsilon}(x,t) \rightarrow \begin{cases} m^{-}_{\beta}, \ \ \ \ \ \ \ \ \ \ \  {\rm in\ }  M_t, \\ m_{\beta}^{+}, \ \ \ \ \ \ \ \ \ \ \  {\rm in\ } \mathbb{R}^3\setminus \overline{M_t},\end{cases}
\end{equation*}
locally uniformly in $\R^N\backslash \partial M_t$,
\end{theorem}

\subsection{Asymptotic expansion near non-characteristic points}

Let us recall that, for all 
 $C^{2}$ function $u: \mathbb{H}^{1} \rightarrow \mathbb{R}$, the Taylor expansion in $\mathbb{H}^{1}$ { (see e.g. \cite{{bo}, {fer}}) is given by
\begin{equation}
 \label{lim}
u(x_0\circ x) = u(x_0) + \big<\nabla_{\mathbb{H}^{1}}u(x_0), \overline{x}\big> + \frac{\partial\,u(x_0)}{\partial x_3}\,  x_3
 + \frac{1}{2} \overline{x}^T \big(D^2_{\mathbb{H}^1} u (x_0)\big)^* \, \overline{x}  + o( \| x \|^{2}_{\mathbb{H}^{1}}),
\end{equation}
where we have used the notation $x=(\overline{x},x_3)\in \R^2\times \R$ for all points in $\mathbb{H}^1$, and the facts that the vector fields are left-translation invariants and  $X_3=\frac{\partial}{\partial x_3}$. Note that by $ \| x \|_{\mathbb{H}^{1}}$ we indicate any norm in $\mathbb{H}^1$  since they are all locally equivalent, e.g. we can take $ \| x \|_{\mathbb{H}^{1}}=\big((x_1^2+x_2^2)^2+16x_3^2\big)^{\frac{1}{4}}$.
As the main difficulty lies in expanding the  convolution,  for the moment we suppress the dependence of functions and surfaces on time, and we set $\beta=1.$ As $\beta$ and $J$ always appear together, this can be done as long as we adjust the $L^1$-norm of $J$ (i.e. the average of $J$).
We apply the previous Taylor expansion to  the function $u(x) = d_{\mathbb{H}^{1}}(x, M)$, and using the simplified following notation: 
$$
 \nu (x_{0})=\nabla_{\mathbb{H}^1} d_{\mathbb{H}^{1}}(x_0, M),
 \;
  l_{3}(x_{0})= \frac{\partial\,}{\partial x_3 }d_{\mathbb{H}^{1}}(x_0, M)
   \; \textrm{and}
 \;
  r= \frac{d_{\mathbb{H}^{1}}(x_{0},  M)}{\varepsilon}, 
$$
and the change of variable $z=(\overline z,z_3)=x_0^{-1}\circ y$, we write for a generic function $m$
 \begin{align*} 
I:=\int_{\mathbb{R}^{3}} J^{\varepsilon}( y^{-1} \circ x_{0}) m \bigg( \frac{d_{\mathbb{H}^{1}}(y,M)}{\varepsilon} \bigg) d y  
 = \int_{\mathbb{R}^{3}} J^{\varepsilon}(z^{-1}) m \bigg( \frac{d_{\mathbb{H}^{1}}(x_0\circ z, M)}{\varepsilon} \bigg) d z=&\\
\!\!\!  \int_{\mathbb{R}^3} \!\!\! J^{\varepsilon}(z^{-1}) m \bigg( r + \nu(x_{0})\cdot \frac{\overline{z}}{\varepsilon} + \varepsilon   l_{3}( x_{0}) \frac{ z_{3}}{\varepsilon^2} + \varepsilon \frac{1}{2} 
\frac{\overline{z}}{\varepsilon} ^{T} \big(D_{\mathbb{H}^{1}}^{2}d_{\mathbb{H}^{1}}(x_{0}, M)
\big)^{*} \frac{\overline{z}}{\varepsilon}  +\frac{o(\|z\|_{\mathbb{H}^{1}}^2)}{\varepsilon}\bigg) dz&,
\end{align*}
where by $\cdot$ we indicate the standard inner product in $\R^2$.

Next we make the change of variables $y= \delta_{\frac{1}{\varepsilon}} (z)$, and use that the determinant of the Jacobian is  $\varepsilon^4$ which cancels the prefactor in  $J^{\varepsilon}(z)=\varepsilon^{-4} J\big( \delta_{\frac{1}{\varepsilon}} (z)
\big)$. Hence, for $y=(\overline y,y_3)\in \R^2\times \R$,  we get
\begin{equation}\label{firstexpansion}
I=
\int_{\mathbb{R}^3}  J(y) m \bigg( r + \nu(x_{0}) \cdot\overline{y} + \varepsilon  l_{3}( x_{0}) y_{3} + \varepsilon \frac{1}{2} \overline{y}^{T} \big(D_{\mathbb{H}^{1}}^{2}d_{\mathbb{H}^{1}}(x_{0}, M)\big)^{*}  \overline{y} +  o(\varepsilon) \|y \|_{\mathbb{H}^{1}}^2 \bigg) dy .
\end{equation}
Here we have also used  the structure $J(y)=J(y_1^2+y_2^2,|y_3|)$ and the fact that in $\mathbb{H}^1$ $y^{-1}=-y$ implies $J(y^{-1})=J(y)$.
Now we expand $m$ to the first order around $r + \nu(x_{0})\cdot \overline{y}$, which implies
\begin{equation}
\label{ful}
\begin{aligned}  
I&=
\int_{\mathbb{R}^3} 
J(y) m \big( r + \nu(x_{0})\cdot\overline{y}  \big) dy +\varepsilon l_{3}( x_{0})  \int_{\mathbb{R}^{3}}  J(y) m' \big( r + \nu(x_{0})\cdot \overline{y}  \big)  y_{3} dy   \\ 
 &+ \frac{\varepsilon }{2} \int_{\mathbb{R}^{3}}  J(y)  m'\big(r + \nu(x_{0})\cdot\overline{y}\big)\;  \overline{y}^{T} \big(D_{\mathbb{H}^{1}}^{2}d_{\mathbb{H}^{1}}(x_{0},M)\big)^{*} \overline{y}\; dy  +  o(\varepsilon),
\end{aligned}
\end{equation}
where we used that $m'$ is bounded and $J$ decays at least exponentially and so it controls all the polynomial terms. 

 By applying Fubini's Theorem in the second integral term in \eqref{ful} to integrate first w.r.t. $y_3$ and using that  $y_3\mapsto J(y)y_3$  is odd, its integral over $\mathbb{R}$ vanishes. This means that the second integral term in \eqref{ful}  is zero.

 We need now to deal with the last integral term above. First, let us note that the only term depending on $y_3$ is the  kernel $J$. By using Fubini's Theorem, we integrate first in $y_3$ and define $\hat{J}(\overline{y})$ as in Definition \ref{kernel}.

 Now we rewrite the remaining 2-dimensional integral in $\overline{y}$ w.r.t. a different basis. 
 It is well known that the Carnot-Carath\'eodory distance solves the eikonal equation (see  Theorem 3.1 of \cite{monti} for $\mathbb{H}^{1}$, and Theorem 2 of \cite{dra}), and so the signed distance $d_{\mathbb{H}^1}$ does.
 This implies that its (non-symmetrized) Hessian at $x_0$  has an eigenvalue zero  w.r.t. the unit length eigenvector  
 $\nu(x_0).$ Recall that, since the Hessian appears only in the form $a^t \big(D_{\mathbb{H}^{1}}^{2}
 d_{\mathbb{H}^{1}}(x_{0}, M)\big)^{*} a$, for some $a\in \mathbb{R}^2$, then  the symmetrized matrix and the non-symmetrized matrix act in the same way.
 \begin{remark}[Change of coordinates]
 We can find a unit vector $\nu^\perp(x_0)$  orthogonal to $\nu(x_0)$ such that $\{ \nu(x_0), \nu^{\perp}(x_0) \}$ is a basis of $\mathbb{R}^2$ with the same orientation as the standard basis, and we denote by $\hat{y}$ the new coordinates in this basis.The other eigenvalue is necessarily the trace, i.e. the mean curvature $k_o(x_0)$ whenever $x_0\in \partial M.$ We denote by $O$ the orthogonal matrix describing this change of basis. 
 This transformation is orthogonal and the determinant of the Jacobian is 1.
 \end{remark}
Since $\hat{J}$ is radially symmetric, so  it is invariant under the change of variables above, which  implies 
 \begin{equation}
 \label{pesci}
 \begin{aligned}
& \int_{\mathbb{R}^{2}}  \hat{J}(\overline{y}) m'\big(r + \nu(x_{0})\overline{y}\big)\;  \overline{y}^{T} \big(D_{\mathbb{H}^{1}}^{2}d_{\mathbb{H}^{1}}(x_{0}, M)\big)^{*} \overline{y}\; d\overline{y}=
 \\
 &\int_{\mathbb{R}^{2}}  \hat{J}(\hat{y}) \ m'(r+ \hat{y}_{1})\, \hat{y}^{T} \, O^{T}\big(D_{\mathbb{H}^{1}}^{2}
d_{\mathbb{H}^{1}}(x_{0}, M)\big)^{*}\, O \,\hat{y}\, d \hat{y}=\\
&\int_{\mathbb{R}^{2}} \hat{J}(\hat{y})\ m'(r+ \hat{y}_{1})\, \hat{y}^{T}\begin{bmatrix} 0 & { 0}\\ 0 & k_{0} (x_0)\end{bmatrix}\hat{y} \;d\hat{y} =\\
& \int_{\mathbb{R}^{2}} \hat{J}(\hat{y}) \ m'(r+\hat{y}_{1})\, k_{0}(x_0)\, \hat{y}_{2}^{2}\, d\hat{y}.
  \end{aligned}
 \end{equation}
 Using \eqref{pesci}  in \eqref{ful} (and the fact that the second integral term vanishes as remarked above), we can write
\begin{equation}
\label{above}
I=\int_{\mathbb{R}^{3}}  J(y) m\big( r + \nu(x_{0})\overline{y}  \big) dy
 +
\frac{\varepsilon}{2}
\int_{\mathbb{R}^{2}} \hat{J}(\hat{y}) m'(r+\hat{y}_{1}) \,k_{0}(x_0) \,\hat{y}_{2}^{2}\, d\hat{y}
+  o(\varepsilon).
\end{equation}
Let us now introduce the following notation:
\begin{equation*}
R(r)=\frac{1}{2}\int_{\mathbb{R}^{3}} \hat{J}(\hat{y}) m'(r+ \hat{y}_{1})   \hat{y}_{2}^{2} d\hat{y},
\quad
\textrm{and}\quad
\overline{J}(\hat y_{1}) = \int_{\mathbb{R}} J(\hat y_{1}^2 + \hat y_{2}^2) d \hat y_2.
\end{equation*}
Then the identity \eqref{above} can be  rewritten as
\begin{equation*}
I= \int_{\mathbb{R}}  \overline{J}(\hat{y}_1)  m ( r + \hat{y}_{1}  ) d\hat{y}_{1}  + \varepsilon k_{0}(x_0) R(r)+ o(\varepsilon).
\end{equation*}
After having expanded the convolution, we apply the previous calculations to the ansatz \eqref{ansatzeucl} adapted to the Heisenberg case inserted into \eqref{Hnonloc}, which leads to
\begin{equation*}
\varepsilon^{-1} \frac{\partial}{\partial\,{t}} d_{\mathbb{H}^{1}}(x,  M) \; m'(r) \! =\! \varepsilon^{-2} \bigg(\! - \widehat m(r) + \tanh \bigg( \! \int_{\mathbb{R}}  \overline{J}(\hat{y}_1) \widehat m ( r + \hat{y}_{1}  ) d\hat{y}_1 + \varepsilon k_{0} R(r) \! \bigg) \bigg) \!+ \! O(1).
\end{equation*}
By expanding the hyperbolic tangent we get
\begin{align*}
\varepsilon^{-1} \frac{\partial}{\partial\,{t}} d_{\mathbb{H}^{1}}(x,  M) \; m'(r)  & = \varepsilon^{-2} \bigg[ -m(r) + \tanh \bigg(  \int_{\mathbb{R}}  \overline{J}(\hat{y}_{1})  m ( r + \hat{y}_{1}  ) d\hat{y}_{1}  \bigg) \\ & + \bigg(1 - \tanh^2 \bigg(  \int_{\mathbb{R}}  \overline{J}(\hat{y}_{1}) m ( r + \hat{y}_{1}  ) d \hat{y}_{1}  \bigg) \bigg) \varepsilon k_{0} R(r)+o(\varepsilon)\bigg] .
\end{align*}
Reasoning as in Section 2, the first non zero term that we obtain is at the level $\varepsilon^{-1}$ and it is
\begin{equation*}
\frac{\partial}{\partial\,{t}} d_{\mathbb{H}^{1}}(x,  M) \, \overline{m}'(r)  = (1-\overline{m}^2) k_{0} R(r) + \mathcal  L(m_{1})(r), 
\end{equation*}
where, see Lemma \ref{linop}, 
\begin{equation*}
{\mathcal  L} f = -f + \big(1-\overline{m}^{2}(r)\big) \int_{\mathbb{R}} \overline{J}(s) f(r+s) ds
\end{equation*}
 is the linearization of non local evolution equation around $\overline{m}$.

As shown in \cite{pres}, and summed up in Lemma \ref{linop}, ${\mathcal L} f =g$ has solution if and only if {$<g,\overline{m}^\prime>_{L^2(\mu)}=0$} where  we recall the notation from Section 2
\begin{equation*}
<f_{1}, f_{2}>_{L^2(\mu)}= \int_{\mathbb{R}^{3}} \frac{1}{1-\overline{m}^2} f_{1}(r) f_{2}(r) dr.
\end{equation*}
Let us split $\hat{R}$ in a perpendicular and parallel component, i.e.
\begin{align*}
\hat{R} & = \left<\hat{R}, \frac{\overline{m}'}{\|\overline{m}'\|_{L^2(\mu)}}\right>_{L^2(\mu)} \frac{\overline{m}'}{\|\overline{m}'\|_{L^2(\mu)}} +  \hat{R} - \left<\hat{R}, \frac{\overline{m}'}{\|\overline{m}'\|_{L^2(\mu)}}\right>_{L^2{(\overline{m})}}\frac{\overline{m}'}{\|\overline{m}'\|_{L^2(\mu)}} \\ &= \hat{R}_{\parallel} + \hat{R}_{\perp}.
\end{align*}
We can choose $m_{1}$ such that ${\mathcal L}(m_{1}) = \hat{R}_{\perp}$ and the parallel term has to cancel with the term multiplying $\overline m ^\prime$ on the left hand side, which gives the evolution law
\begin{equation*}
V_0=\frac{\partial}{\partial_{t}} d_{\mathbb{H}} = k_{0} \theta,
\end{equation*}
where 
\begin{eqnarray}\label{theta}\theta  &=&\left<\hat{R}, \frac{\overline{m}'}{\|\overline{m}'\|_{L^2(\mu)}}\right>_{L^2{(\mu)}}\!\!\!\!  \!\!\!=\frac{\beta}{N} \int_{\mathbb{R}^{3}} \overline{m}'(r) \overline{m}'(r+r_{1}) \hat{J}(r_{1}^2 + s^2) s^2 dr dr_{1} ds\\ \nonumber N &=&\int_{- \infty}^{\infty} \frac{(\overline{m}'(r))^2}{1- \overline{m}^2(r)} dr.
\end{eqnarray}

\subsection{Characteristic points and numerical results}

At characteristic points, the horizontal mean curvature is not defined due to the absence of the horizontal normal at these points. On the other hand, the nonlocal evolution equation \eqref{Hnonloc} is defined even at characteristic points. In this section we study a version of \eqref{Hnonloc} numerically, starting from a level set function for a sphere in the homogeneous norm, which has characteristic points in the north and south pole. Let us sum up the findings.

The balls shrinks everywhere, even at the characteristic points,  with continuous but non-uniform velocity. The profile of $m^\varepsilon$ near the characteristic points is steeper than away from them, because the scaling in the direction of the $x_3-$axis is $\varepsilon^2,$ not $\varepsilon.$
We now explain the algorithm and the results in more detail.

As kernel we take the heat kernel of the Heisenberg group. It is not compactly supported, but due to its  exponential decay this is a minor issue. As the convolution in the Heisenberg group is not the standard Euclidean convolution, we cannot compute it by fast Fourier transform. The choice of the heat kernel allows us to solve the heat equation on the Heisenberg group with linear finite elements instead. Then we solve the ODE \eqref{Hnonloc} by a semi-implicit Euler scheme: The linear part is implicit, the nonlinear (tanh) part explicit.
This means we perform the following steps for a time step size $\Delta t$ and $k\in \mathbb{N},$ $\varepsilon>0$ with $\Delta t\varepsilon^{-2}<1,$
$\beta>1:$

\begin{enumerate}
\item
Convolution with heat kernel by solving heat equation:
\begin{eqnarray*}
\frac{\partial}{\partial_t} v(t,x)&=&\Delta_{{\mathbb H}^1}v(t,x)\\
v(0,x)&=&m((k-1)\Delta t,x)
\end{eqnarray*}
\item
$$\begin{aligned}
m(k\Delta t,x)=\frac{1}{1+\varepsilon^{-2}\Delta t}m((k-1)\Delta t,x)+\frac{\varepsilon^{-2}\Delta t}{1+\varepsilon^{-2}\Delta t}
\tanh(\beta v(\varepsilon^2,x)) .
\end{aligned}
$$
\end{enumerate}
Note that the equation in step 2 can be reformulated: take $0<\delta=\frac{\varepsilon^{-2}\Delta t}{1+\varepsilon^{-2}\Delta t}<1,$ then
$$
m(k\Delta t,x)=(1-\delta)m((k-1)\Delta t,x)+\delta
\tanh(\beta v(\varepsilon^2,x)).
$$This indicates that  the limiting case $\delta\to 1,$ $\beta\to \infty$ is related to a diffusion-concentration algorithm, as, for example, the  Bence–Merriman–Osher algorithm \\ (Merriman et al., 1992 [42]).  A sub-Riemannian analogue was introduced and the convergence at the non-characteristic case discussed in  \cite{cap}. On the other hand, for $\beta$ fixed and $\delta\to 0,$ the algorithm converges to the solution of the nonlocal equation, i.e. the model by Ermentrout and Cowan \cite{co}, so our algorithm interpolates between the two.

Ferrari, Manfredi and Liu \cite{fer} have shown 
 comparison principles for the mean curvature flow in $\mathbb{H}^1$ in the case of  compact surfaces if at least one of the surfaces is rotationally symmetric around the $x_3$-axis, even  in the presence of characteristic points. This allows to deduce that any compact surface  shrinks to a point in finite time. 
Recalling the homogeneous distance in $\mathbb{H}^1$ is given by
\begin{equation*}
\|x\|_{\mathbb{H}^1}=\left( (x_{1}^2+x_{2}^2)^2+16x_{3}^2\right)^{\frac{1}{4}},
\end{equation*}
in \cite{fer} the authors give an explicit solution starting from the hypersarface 
$$\partial B_r:=\{x\in \mathbb{H}^1:\ \|x\|_{\mathbb{H}^1}=r\},$$ which has characteristic points at its intersection with the $x_3$-axis.
Unlike in the Euclidean case, this solution is not self-similar, but it is symmetric around the $x_3$-axis and thus unique. We test our algorithm against this exact solution. 
 The choice $\beta=1.2$ leads to $m_\beta=0.6585.$  Then we find the profile $\overline m$ (which depends on $\beta$) by choosing an initial condition with planar symmetry and letting the evolution run until the solution stabilises.

The next step is to determine the factor $\theta.$  
A cylinder revolving around the $x_3$-axis evolves according to the 2-dimensional Euclidean mean curvature flow. Starting from such a geometry, our code produces results consistent with a radius $r(t)$ evolving like
$$
\dot r(t)=-\frac{\theta}{r(t)}\quad \Longrightarrow \quad \frac{1}{2}\left(r^2(t)\right)= \frac{1}{2}\left(r^2(0)\right)-{\theta}t.
$$
For our choice of $\beta$ and the chosen kernel, we find 
$\theta=0.56561$ by a linear regression, see Figure \ref{fig:linreg}.

\begin{figure}[h!] \label{fig:linreg} \begin{center}
\includegraphics[scale=0.175]{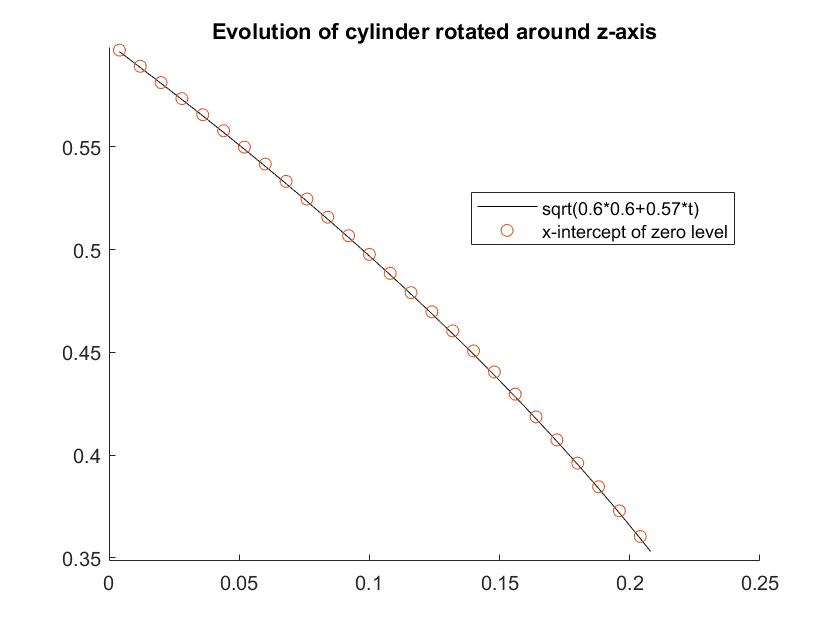}\quad \includegraphics[scale=0.175]{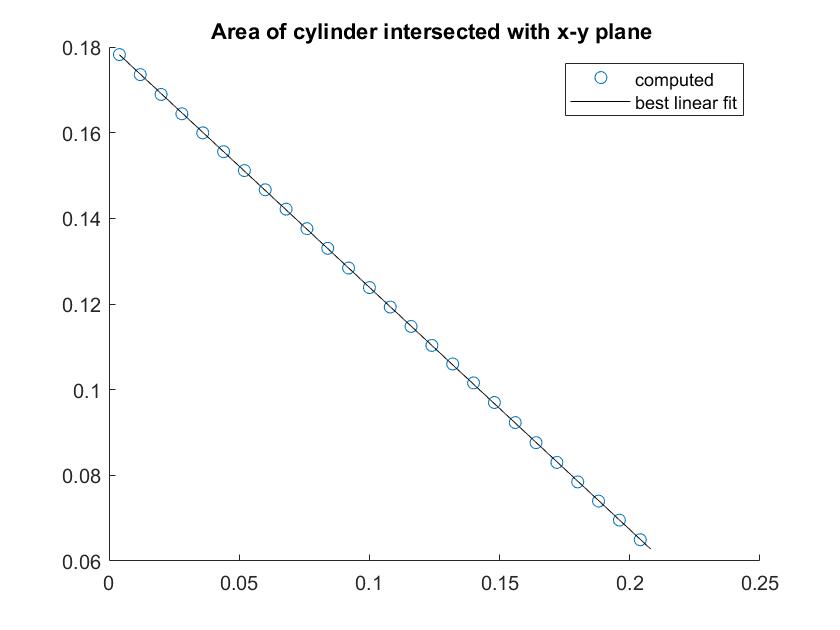}
\end{center}
\caption{Computation of $\theta$, $x$-axis time, $y$-axis radius (left) and area (right) of circle}
\end{figure}

Then we compute the nonlocal evolution equation with the heat kernel and initial condition 
$\overline m((\|x\|_{\mathbb{H}^1}-r)/\varepsilon),$ which has $\partial B_r$ as its zero level set.

The exact solution at time $t$ is obtained as (see \cite{fer})
$$
\partial M_t=\{x\in \mathbb{H}^1:\ (x_{1}^2+x_{2}^2)^2+12\theta t (x_{1}^2+x_{2}^2)+16 x_{3}^2+12 (\theta t)^2=r^4\},
$$
simply by a time-rescaling  in the formula in \cite{fer}.
The parameters are $\beta=1.2,$ $\varepsilon=0.1$ and $r=1.2.$ In Figure \ref{fig:Homball} we show a good qualitative agreement and the inward motion of the characteristic points.

\begin{figure}[h!]\label{fig:Homball}\begin{center}
\includegraphics[scale=0.3]{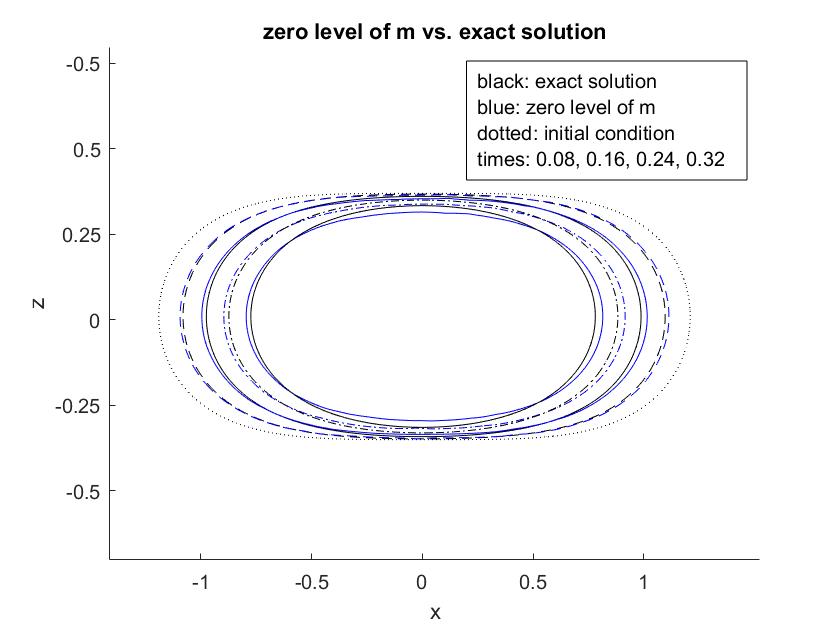}
\end{center}
\caption{ Zero level sets of $m^\varepsilon(t, x_{1},0,x_{3})$ for different times, both for exact solution of the Heisenberg MCF (black) and the level set of the numerically solved nonlocal evolution equation (blue).}
\end{figure} 

In  Figure \ref{fig:profile} we can note a much steeper profile at the  characteristic points, as the kernel there effectively scales with $\varepsilon^2.$ The profiles have both been translated to be centred at zero.

\begin{figure}[h!]\label{fig:profile}\begin{center}
\includegraphics[scale=0.25]{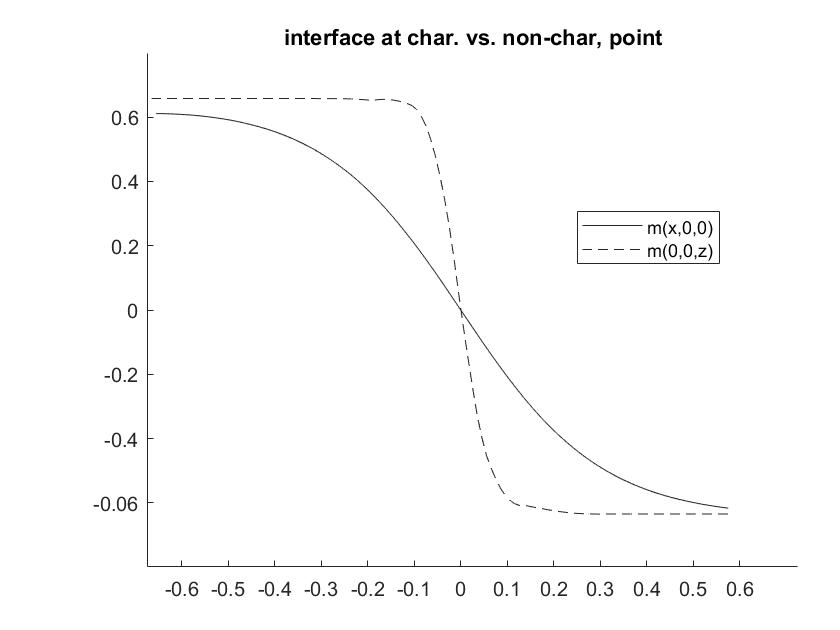}
\end{center}
\caption{Graph of $m^\varepsilon(t, x_{1},0,0)$ (solid) and $m^\varepsilon(t, 0,0,x_{3})$ (dashed)  at $t=0.32$}
\end{figure} 

These computations indicate that the zero level set of the solution to the nonlocal evolution equation approximates the (horizontal)  mena curvature flow even near  the characteristic points. For practical applications, one would not use the heat kernel, but instead choose a short range kernel in order to improve speed significantly, and one would also consider parallelization, which is convenient for a kernel with relatively short range.

\section{Extension to the $SE(2)$ group}

In this section we show how to extend the result to the group $SE(2)$. Since the Heisenberg group is the tangent cone to the roto-translation group, they present the same local approximation. Here we briefly introduce the group and refer to \cite{cisar2} for more results. We will identify the group with the space $\mathbb{R}^2 \times S^1$, and denote by $x=(x_1,x_2, \theta)$ the general point of the space, where $ (x_1, x_2) \in \mathbb{R}^2$ and $\theta \in S^1$ In the future we write  $\hat x=(x_1,x_2)\in \mathbb{R}^2.$ We will also consider the choice of vector fields 
$$
Y_1(x) =
\left(\begin{array}{c}\cos\theta\\
\sin\theta\\0\end{array}\right),
\quad 
Y_2(x)= 
\left(\begin{array}{c}0\\
0\\1\end{array}\right),$$
and $Y_3=[Y_1,Y_2]$.
A natural composition law is defined by the composition of matrices
and the vector fields are left-invariant w.r.t. this non-commutative group law. 

Note also that 
$$
Y_3(x) = [Y_1, Y_2](x) $$
is linearly independent of the other two vector fields,
hence the H\"ormander condition is satisfied also in this setting. Consequently a distance function can be defined as in \eqref{dista}.

The fact that the Heisenberg group is the tangent cone to $SE(2)$ is explained in Chapter 9.4.3 of \cite{LeDonne}, e.g. Theorem 9.4.6.

In practice, we do  not have a dilation in the space, but we can define a dilation on the Lie algebra 
$$
\delta_\lambda (a) = (\lambda a_1, \lambda a_2, \lambda^2 a_3).
$$
Here and in the sequel we will denote $a$ the coordinates on the Lie algebra, and with $x$ the elements of the 
group. We use the exponential map to define local dilation around a single point in the group. Precisely we call canonical coordinates around a point $x$: 
$$
\Theta_x(y) = a \text{ iff } y = \exp(a_1 Y_1 + a_2 Y_2 + a_3 Y_3)(x),
$$
 and consequently we can locally define dilation around a point $x$: 
$$
\delta_{\lambda, x}(y) = \Theta_x \delta_	\lambda \Theta^{-1}_x (y), \quad \textrm{for}\; \lambda\geq 0.
$$

Assume that $J:\mathbb{R}^3 \to \mathbb{R}$ is compactly supported. We define 
$$J^\varepsilon(a) = \frac{1}{\varepsilon^4}J(\delta_{\varepsilon}(a)).$$

The Taylor expansion of a general function $u$ of class $C^2$ with respect to the present vector fields 
is the exact analogous to the one recalled in \eqref{lim}: 
$$ u(\exp( a_i X_i) (x_0)) = u(x_0) 
+ \sum_{i=1}^3 a_i X_i u(x_0) + \frac{1}{2}
a^T(D_S^2u(x_0) )^*
x_0 a + o(||a||^2),$$
where $D_S^2$ is the horizontal Hessian defined exactly as in \eqref{lim} but with respect to the present vector fields, and $||\quad ||$ is any homogenous norm.
Arguing as in \eqref{firstexpansion}, we can write for a generic function $m$

\begin{align*}I :&=
\int_{\mathbb{R}^3}J^\varepsilon(a )m\Big(\frac{d(\exp(a)(x_0),M)}{\varepsilon}
\Big)da\\
 &=
\int_{\mathbb{R}^{3}} 
J^\varepsilon(a)
m\Big( \frac{d( x_0,M)}{\varepsilon}  
+ \sum_{i=1}^3 \frac{a_i}{\varepsilon} X_i u(x_0) + 
\frac{1}{2}
\frac{a^T}{\varepsilon}(D^2u(x_0) )^*
 a + o \bigg(\frac{||a||^2}{\varepsilon}\bigg)\bigg)da,
 \end{align*}
and expand with the change of variable 
$b = \delta_{1/\varepsilon}(a)$. Note that we are on the algebra, hence this is the same dilation as in Heisenberg, giving
$$I=\!\!
\int_{\R^3}\!\!
J(b)
m\Big( \frac{d( x_0,M)}{\varepsilon}  
+ \sum_{i=1}^2 b_i X_i u(x_0) + 
\varepsilon b_3 X_3 u(x_0)+
\frac{\varepsilon}{2} b^T(D^2u(x_0) )^*
b + o(\varepsilon ||b^2||)\Big)db. $$

We expand $m$, obtaining
\begin{align*}
I&=
\int_{\R^3} 
J(b)
m\Big( \frac{d( x_0,M)}{\varepsilon}  
+ \sum_{i=1}^2 b_i X_i u(x_0)\Big) db\\
 &+ \varepsilon X_3 u(x_0) \int_{\R^3} J(b)
 b_3 m'\Big( \frac{d( x_0,M)}{\varepsilon}  
+ \sum_{i=1}^2 b_i X_i u(x_0)\Big) db \\
&+
\frac{\varepsilon}{2} \int_{\R^3}   J(b)
 m'\Big( \frac{d( x_0,M)}{\varepsilon}  
+ \sum_{i=1}^2 b_i X_i u(x_0)\Big)  
 b^T(D^2u(x_0) )^*
b\, db + o(\varepsilon ).
\end{align*}

Since $J$ is defined on the algebra, which coincides with $\mathbb{R}^3$, exactly the same argument works in this setting.

Since $ J(b)b_3$ is odd, its integral over $\mathbb{R}$ vanishes,
that means that the second integral term above is zero.

\begin{align*}
I&=
\int_{\R^3} 
J(b)
m\Big( \frac{d( x_0,M)}{\varepsilon}  
+ \sum_{i=1}^2 b_i X_i u(x_0)\Big) \\
&+
\frac{\varepsilon}{2} \int_{\R^3}    J(b)
 m'\Big( \frac{d( x_0,M)}{\varepsilon}  
+ \sum_{i=1}^2 b_i X_i u(x_0)\Big)  
 b^T(D^2u(x_0) )^*
b\, db + o(\varepsilon ). 
\end{align*}

We need now to deal with the last integral term above. First, let us note
that the only term depending on $b_3$ is the kernel $J$. We integrate with respect to $b_3$ and reduce to $\hat J:\R^2\to \R$ defined as 
$$\hat J(\hat b) = \int_\R J(b) db_3, \quad \textrm{for} \; \hat b=(b_1,b_2)\in \R^2.$$
We obtain 
\begin{align*}
  I &  = \int_{\R^3} 
 J( b)
m\Big( \frac{d( x_0,M)}{\varepsilon}  
+ \sum_{i=1}^2 b_i X_i u(x_0)\Big) d b \\
&+
\frac{\varepsilon}{2} \int_{\R^2}  \hat  J(\hat b)
 m'\Big( \frac{d( x_0,M)}{\varepsilon}  
+ \sum_{i=1}^2 b_i X_i u(x_0)\Big)  
 \hat b^T(D^2u(x_0) )^*
\hat b \,d\hat b + o(\varepsilon ). 
\end{align*}

Arguing as in the Heisenberg case we reach 
the following equation
\begin{align*}
 I & = \int_{\R^3} 
J(b)
m\Big( \frac{d( x_0,M)}{\varepsilon}  
+ \sum_{i=1}^2 b_i X_i u(x_0)\Big) db \\
&+\frac{k(x_0) \varepsilon}{2} \int_{\R^2}  \hat  J(\hat b)
 m'\Big( \frac{d( x_0,M)}{\varepsilon}  
+ b_1 \Big)  
 b_2^2 d\hat b + o(\varepsilon ),
 \end{align*}
which is the exact analogous of Equation \eqref{above}.
From this point we can
 proceed as in the Heisenberg case in order to conclude.

\end{document}